\documentclass[11 pt]{amsart}
\usepackage{amsmath}
\usepackage[T1]{fontenc}
\usepackage{textcomp} 
\DeclareEncodingSubset{TS1}{hlh}{1}  
\usepackage{mathalfa} 
\usepackage{amsmath, amsthm}
\usepackage{amssymb}
 \usepackage{amsaddr}
 \usepackage{graphicx}
\usepackage[utf8]{inputenc}
\usepackage{rotating}
\usepackage{tikz}
\usepackage{rotfloat}
\usepackage{hyperref}
\usepackage{setspace}
\theoremstyle{definition}

\theoremstyle{remark}

\numberwithin{equation}{section}

\usepackage{lscape}
\usepackage{setspace}
	
\usepackage{lscape}

\begin{document}

\title{Rogers-Ramanujan continued fraction and approximations to $\mathbf{2\pi}$}
\author{Rajeev Kohli}
\date{\today}
\thanks{{\it Address:}  972 Kravis Hall, 
Columbia University, 
665 W 130th St, New York, NY 10027}

\maketitle

\begin{abstract}
We observe that certain famous evaluations of the Rogers-Ramanujan continued fraction $R(q)$ are close to $2\pi-6$ and $(2\pi-6)/2\pi$, and that   
$2\pi-6$ can be expressed by a Rogers-Ramanujan continued fraction in which $q$ is very nearly equal to $R^5(e^{-2\pi})$. The value of $-{5\over \alpha}\ln R(e^{-2\alpha \pi})$ converges to $2\pi$ as $\alpha$ increases. For $\alpha=5^n$, a modular equation by Ramanujan provides recursive closed-form expressions that approximate the value of $2\pi$, the number of correct digits increasing by a factor of five each time $n$ increases by one. If we forgo closed-form expressions, a modular equation by Rogers allows numerical iterations that converge still faster to $2\pi$, each iteration increasing the number of correct digits by a multiple of eleven.
\end{abstract}


\section{Introduction}
\noindent
Ramanujan evaluated the Rogers-Ramanujan continued fraction $R(q)$ for several different values of $q$. We show that two of his evaluations are close to $2\pi-6$ and $(2\pi-6)/2\pi$, and can be combined with a third to obtain another close approximation for $2\pi-6$. We use these relations to express $2\pi$ by a Rogers-Ramanujan continued fraction $R(q)$ when $q=\rho^5 R^5(e^{-2\pi})$, where $\rho$ is nearly equal to one. We use a modular equation by Ramanujan to obtain closed-form approximations for $2\pi$ in the form $-\ln R(e^{2\cdot 5^n \pi})/5^{n-1}$, where $n\ge 0$. Increasing the value of $n$ by one increases the number of correct digits of $2\pi$ by a factor of five. The closed form has the benefit that it avoids the loss in accuracy associated with iterated numerical computation. But if we forgo the closed form and use numerical iterations, then a modular equation due to Rogers can be used to obtain an eleven-fold increase in the number of digits of $2\pi$ per iteration.

\section{$R(q)$ and approximations to $\pi$}

\bigskip\noindent
We investigate approximate relationships between $2\pi$ and the Rogers-Ramanujan continued fraction 

$$R(q)={q^{1/5}\over 1+\ }{q\over 1+\ }{q^2\over 1+\ }{q^3\over 1+\ }\dots,\ \text{where}\ |q|<1.$$

\medskip\noindent
Ramanujan recorded several values of $R(q)$ in his notebooks. His first two letters to  Hardy (\cite{Hardy}) included the following three evaluations:

$$R(e^{-2\pi})=\sqrt{\phi^2+1}-\phi$$

$$R^5(e^{-2\pi/\sqrt{5}})=\sqrt{\phi^{10}+1}-\phi^5$$

and

$$R\left(e^{-2\pi\sqrt{5}}\right)
={1-\phi R\left(e^{-2\pi/\sqrt{5}}\right)\over \phi+R\left(e^{-2\pi/\sqrt{5}}\right)}
={1-\phi \sqrt[5]{{\sqrt{\phi^{10}+1}-\phi^5}}\over \phi +\sqrt[5]{{\sqrt{\phi^{10}+1}-\phi^5}}}$$

\medskip\noindent
where $\phi=(\sqrt{5}+1)/2$. We make the following observations.

\vskip 0.25in\noindent
(1) $R(e^{-2\pi})$ is very close to $2\pi-6$:

$$R(e^{-2\pi})=2\pi-6+8.937\dots \times 10^{-4}\ .\quad \eqno(1)$$

\vskip 0.25in\noindent
(2) $R^5(e^{-2\pi/\sqrt{5}})$ is still closer to $(2\pi-6)/2\pi$:

$$R^5(e^{-2\pi/\sqrt{5}}) = {2\pi-6\over 2\pi} + 7.6641082\dots \times 10^{-5}\ . \eqno(2)$$

\medskip\noindent
Equation (2) implies that

$${6R^5(e^{-2\pi/\sqrt{5}})\over {1-R^5(e^{-2\pi/\sqrt{5}})}}
=2\pi-6-5.042376378\dots \times 10^{-4}\ . \eqno(3)$$

\vskip 0.25in\noindent
(3) $R(e^{-2\pi\sqrt{5}})$ is nearly equal to $(2\pi-6)\sqrt{{2\pi -6\over 2 \pi}}$\ : 

$$R(e^{-2\pi\sqrt{5}})=(2\pi-6)\sqrt{{2\pi -6\over 2 \pi}}+8.97985\dots \times 10^{-5}\ . \eqno(4)$$

\medskip\noindent
Equation (4) implies that

$${R(e^{-2\pi\sqrt{5}})\over \sqrt{R^5(e^{-2\pi/\sqrt{5}})}}=2\pi-6+6.6442631\dots \times 10^{-4}\ . \eqno(5)$$

\vskip 0.25in\noindent
Note that

$${6R^5(e^{-2\pi/\sqrt{5}})\over {1-R^5(e^{-2\pi/\sqrt{5}})}}
<2\pi-6
<{R(e^{-2\pi\sqrt{5}})\over \sqrt{R^5(e^{-2\pi/\sqrt{5}})}}
<R(e^{-2\pi})
<e^{-2\pi/5}\ .\ \eqno(6)$$

\medskip\noindent
Table 1 shows the value of each term, and the differences in the values of successive  terms,  in equation (6). 

\begin{table}[]
\caption{Values of five terms and their differences}
\label{values}
\begin{tabular}{lll}
\hline 
&&Difference between\\
&Value&successive terms\\
Term&(10 digits)&(10 digits)\\

\hline \\

\bigskip
$x_1={6R^5(e^{-2\pi/\sqrt{5}})\over {1-R^5(e^{-2\pi/\sqrt{5}})}}$ &0.2826810695 & \\

\bigskip
$x_2=2\pi-6$ & 0.2831853072 &$x_2-x_1=0.0005042376$\\

\bigskip
$x_3={R(e^{-2\pi\sqrt{5}})\over \sqrt{R^5(e^{-2\pi/\sqrt{5}})}}$&0.2838497335  &$x_3-x_2=0.0006644263$  \\

\bigskip
$x_4=R(e^{-2\pi})$ &0.2840790438  &$x_4-x_3=0.0002293103$    \\

\bigskip
$x_5=e^{-2\pi/5}$ &0.2846095433 &$x_5-x_4=0.0005304994$  \\

\hline
\end{tabular}
\end{table}

\vskip 0.25in\noindent
The relation $e^{-2\pi/5}\approx R(e^{-2\pi})\approx 2\pi-6$ implies that

$$R(q)=2\pi-6\eqno(7)$$

\medskip\noindent
when

$$q^{1/5}=\rho R(e^{-2\pi})=\rho \left(\sqrt{\phi^2+1}-\phi\right),$$

\medskip\noindent
where

$$\rho=0.9997370833\dots \approx 1.$$

\medskip\noindent
Equation (7) ``inverts'' $R(e^{-2\pi})$, expressing $2\pi$ by $6+R(q)$ when $q$ is proportional (and very nearly equal) to $R^5(e^{-2\pi})$:

$$2\pi=6+{\rho R(e^{-2\pi})\over 1+\ }\ {(\rho R(e^{-2\pi}))^5\over 1+\ }\ {(\rho R(e^{-2\pi}))^{10}\over 1+\ }\ {(\rho R(e^{-2\pi}))^{15}\over 1+\ }\dots.$$

\medskip\noindent
Similarly, equations (3) and (5) can be used to express $2\pi$ by Rogers-Ramanujan continued fractions in terms of $R^5(e^{-2\pi/\sqrt{5}})$ and $R(e^{-2\pi\sqrt{5}})$.

\vskip 0.25in\noindent
One interpretation of equations (1), (3) and (5) is that $6+R(e^{-2\pi})$, $6/(1-R^5(e^{-2\pi/\sqrt{5}}))$ and $6+ (R(e^{-2\pi\sqrt{5}})/{R^{5/2}(e^{-2\pi/\sqrt{5}})})$ are the perimeters of ellipses, each very nearly a unit circle, with minor axes equal to one and major axes slightly greater than one. In each case, we can estimate the major axis using Ramanujan's approximation formula for the perimeter of an ellipse \cite{Ramanujan1}

$$p\approx \pi (a+b)\left(1+{3\lambda^2\over 10+\sqrt{4-3\lambda^2}}\right),$$

\medskip\noindent
where $a=1+d$ is the major axis, $b=1$ the minor axis and 

$$\lambda={a-b\over a+b}={d\over 2+d}\ .$$

\medskip\noindent
For example, if the perimeter is $6+R(e^{-2\pi})$, then 

$$d\approx 0.0002844725721532\dots$$

\medskip\noindent
which is very nearly equal to 

$${R(e^{-2\pi})\over 1000}=0.0002840790438404\dots$$

\vskip 0.25in\noindent
Let $q=e^{-2\pi\alpha}$, where $\alpha\ge 0$. Then 

$$\lim_{\alpha\rightarrow \infty} {R(e^{-2\pi\alpha})\over e^{-2\pi\alpha/5}}=\lim_{\alpha\rightarrow \infty}  {1\over 1+\ }{e^{-2\pi\alpha}\over 1+\ }{e^{-4\pi\alpha}\over 1+\ }{e^{-6\pi\alpha}\over 1+\ }\dots=1$$

\medskip\noindent
and

$$2\pi \approx -{5\over \alpha}\ln R(e^{-2\alpha\pi}), \eqno(8)$$

\medskip\noindent
where the approximation improves as $\alpha$ increases. 

\bigskip\noindent
Ramanujan \cite{Ramanujan1} obtained several approximations to $\pi$ in terms of the logarithms of surds using expressions for $G_n$ and $g_n$. Equation (8) provides approximations to $2\pi$ directly in terms of the logarithm of $R(q)$. Let $\alpha=5^n$, where $n\ge 0$. Then we can recursively obtain closed-form expressions for $R(e^{-2\cdot 5^n \pi})$. Let $u=R(q)$ and $v=R(q^5)$. Ramanujan obtained the following relation in his lost notebook (\cite{AB}, Entry 3.2.13, p. 93):

$$u^5=v{1-2v+4v^2-3v^3+v^4\over 1+3v+4v^2+2v^3+v^4}\ .\eqno(9)$$

\medskip\noindent
Chan, Cooper and Liaw \cite{CCL} used equation (9) to obtain

$$v={1-\phi Y\over \phi + Y}, \eqno(10)$$

\medskip\noindent
where

$$Y=\sqrt[5]{{1-\phi^5u^5\over \phi^5+u^5}}.\eqno(11)$$

\bigskip\noindent
Let $n=0$. Then $\alpha=5^n=1$ and $R(e^{-2\alpha\pi})=R(e^{-2\pi})$. Equation (8) gives

$$2\pi \approx -{5\over 1} \ln R(e^{-2\pi}),$$

\medskip\noindent
which has the approximation error

$$2\pi+5 \ln \left(\sqrt{\phi^2+1}-\phi\right)=-9.3284736\dots \times 10^{-3}.$$

\bigskip\noindent
Let $n=1$. Substituting $u={R(e^{-2\pi})}=\sqrt{\phi^2+1}-\phi$ into equation (11) gives

$$Y=\sqrt[5]{{1-\phi^5(\sqrt{\phi^2+1}-\phi)^5\over \phi^5+(\sqrt{\phi^2+1}-\phi)^5}}=\sqrt[5]{3\left(\sqrt{\phi^2+1}-1\right)-\phi^2}.$$

\bigskip\noindent
Substituting for $Y$ into equation (10) yields

$$v=R(e^{-10\pi})={1-\phi Y\over \phi+Y}\qquad\qquad\qquad\ $$

$$=\frac{1-\phi \sqrt[5]{3\left(\sqrt{\phi^2+1}-1\right)-\phi^2}}{\phi +\sqrt[5]{3\left(\sqrt{\phi^2+1}-1\right)-\phi^2}}\ .\eqno(12)$$

\medskip\noindent
It follows from equation (8) that 

$$2\pi \approx -{5\over 5} \ln R(e^{-10\pi}),$$

\medskip\noindent
which has the approximation error 

$$2\pi+\ln R(e^{-10\pi})=-2.2711010\dots \times 10^{-14}\ .$$

\bigskip\noindent
Let $n=2$. Substituting $u=R(e^{-2\cdot 5\pi})$ into equations (11) and (10) gives

$$Y=\sqrt[5]{{(\phi + X)^5-\phi^5(1-\phi X)^5\over \phi^5 (\phi+X)^5 + (1-\phi X)^5}}={A\over B}$$

\medskip\noindent
and

$$v=R(e^{-50\pi})={1-\phi Y\over \phi+Y}={B-\phi A\over A+\phi B},$$

\medskip\noindent
where

$$X=\sqrt[5]{3\left(\sqrt{\phi^2+1}-1\right)-\phi^2},$$

\medskip\noindent
$$A=\sqrt[5]{{\left(\phi +\sqrt[5]{3\left(\sqrt{\phi^2+1}-1\right)-\phi^2
}\right)^5-\phi ^5 \left(1-\phi \sqrt[5]{3\left(\sqrt{\phi^2+1}-1\right)-\phi^2} \right)^5}}$$

\medskip\noindent
and

\medskip\noindent
$$B=\sqrt[5]{{\phi ^5 \left(\phi +\sqrt[5]{3\left(\sqrt{\phi^2+1}-1\right)-\phi^2}\right)^5+\left(1- \phi \sqrt[5]{3\left(\sqrt{\phi^2+1}-1\right)-\phi^2}\right)^5}}\ .$$

\bigskip\noindent
Using equation (8) gives

$$2\pi \approx -{5\over 25}\ln R(e^{-50\pi}),$$

\medskip\noindent
which has the approximation error

$$2\pi + {1\over 5}\ln R(e^{-50\pi})=-1.20840441\times 10^{-69}.$$

\medskip\noindent
Table 2 shows the number of leading digits of $2\pi$ obtained for $0\le n\le 6$. The number of correct digits increases by a factor of five each time $n$ increases by one. Thus, $n=8$ yields over one million digits of $\pi$, and $n=23$ yields more than 3.5 times the 100 trillion digits of $\pi$ calculated by Google \cite{Iwao} using Chudnovsky's algorithm \cite{Chudnovsky}, \cite{Milla}. Although it is too long to be written down, the closed form has the benefit that it avoids the loss in accuracy associated with iterated numerical computation. Other methods for calculating the digits of $\pi$ include Chan et al.'s \cite{CCL} quintic approximation for $1/\pi$; the series for $1/\pi$ obtained by the 
Borweins \cite{BB3}, which yields 50 digits of $\pi$ per term; and Berndt and Chan's series \cite{Berndt}, which obtains about 73 or 74 digits of $\pi$ per term.

\bigskip
\begin{table}[htp]
\caption{Number of correct digits $k$ in the approximation to $2\pi$ as a function of $n$.}
\begin{center}
\begin{tabular}{ll}
\hline
$n$&$k$\\
\hline
0&$3$\\
1&$14$\\
2&$69$\\
3&$342$\\
4&$1,706$\\
5&$8,528$\\
6&$42,637$\\
\hline
\end{tabular}
\end{center}
\label{$k$}
\end{table}%

\bigskip\noindent
If we forgo closed-form expressions, we can use modular relations to obtain greater approximation accuracy in each step of a numerical iteration. Let $u=R(q)$ and $v=R(q^{11})$. The following modular relation by Rogers \cite{Rogers} is a consequence of two modular equations obtained by Ramanujan \cite{BCHKSS}:

$$uv\left(1-11u^5-u^{10}\right)\left(1-11v^5-v^{10}\right)=\left(u-v\right)^{12}.\eqno(13)$$

\medskip\noindent
To illustrate the iteration, we begin with $u=R(e^{-10\pi})$, which has the value in equation (12). We substitute for $u$ and numerically solve for $v=R(e^{-110\pi})$ in equation (13). From equation (8),

$$2\pi \approx -{5\over 55}\ln R(e^{-110\pi}),$$

\medskip\noindent
which has the approximation error 

$$-{1\over 11}\ln R(e^{-110\pi})-2\pi=7.5371714126\dots \times 10^{-152}\ .$$

\medskip\noindent
In the second iteration, we substitute $u=R(e^{-110\pi})$ and numerically solve for $v=R(e^{-11\times 110\pi})$ in equation (13). Then we use equation (8) to obtain 

$$2\pi \approx -{5\over 605}\ln R(e^{-1210\pi}),$$

\medskip\noindent
which has the approximation error 

$$-{1\over 121}\ln R(e^{-1210\pi})-2\pi=1.0515416546\dots \times 10^{-1653}\ .$$

\medskip\noindent
Thus, the first iteration gives 152 correct digits, and the second 1653 correct digits, of $2\pi$. Each additional iteration multiplies the number of correct digits by a factor of eleven.

\end{document}